\documentclass[%
 aip,
 amsmath,amssymb,
 reprint,%
]{revtex4-1}

\usepackage{graphicx}% Include figure files
\usepackage{dcolumn}% Align table columns on decimal point
\usepackage{bm}% bold math
\usepackage[utf8]{inputenc}
\usepackage[T1]{fontenc}
\usepackage{mathptmx}
\usepackage{etoolbox}

\usepackage{graphicx}% Include figure files
\usepackage{dcolumn}% Align table columns on decimal point
\usepackage{bm}% bold math
\usepackage{overpic}
\usepackage{tabularx}
\usepackage{multirow}
\usepackage{makecell}

\usepackage{amsmath}
\DeclareMathOperator{\NN}{NN}
\DeclareMathOperator{\ODESolve}{ODESolve}
\DeclareMathOperator{\vectorize}{vec}

\begin{document}

\preprint{APS/123-QED}

\title{Structural Inference of Networked Dynamical Systems with Universal Differential Equations} % This could be anything related to topology, sys. ID., etc

\author{J. Koch}
 \email{james.koch@pnnl.gov}
\author{Z. Chen}
\author{A. Tuor}
\author{J. Drgona}
\author{D. Vrabie}
\affiliation{ 
Pacific Northwest National Laboratory, Richland, WA
}%

\date{\today}% It is always \today, today,
             %  but any date may be explicitly specified

\begin{abstract}

Networked dynamical systems are common throughout science in engineering; e.g., biological networks, reaction networks, power systems, and the like. For many such systems, nonlinearity drives populations of identical (or near-identical) units to exhibit a wide range of nontrivial behaviors, such as the emergence of coherent structures (e.g., waves and patterns) or otherwise notable dynamics (e.g., synchrony and chaos). In this work, we seek to infer (i) the intrinsic physics of a base unit of a population, (ii) the underlying graphical structure shared between units, and (iii) the coupling physics of a given networked dynamical system given observations of nodal states. These tasks are formulated around the notion of the \textit{Universal Differential Equation}, whereby unknown dynamical systems can be approximated with neural networks, mathematical terms known \textit{a priori} (albeit with unknown parameterizations), or combinations of the two. We demonstrate the value of these inference tasks by investigating not only future state predictions but also the inference of system behavior on varied network topologies. The effectiveness and utility of these methods is shown with their application to canonical networked nonlinear coupled oscillators.
\end{abstract}

\maketitle

\begin{quotation}
Data-driven modeling of dynamical systems has experienced a surge in method development in concert with advances in machine learning and artificial intelligence. In this work, we restrict our scope specifically to the problem of data-driven modeling of networked systems; a challenging problem in which one needs to elicit not only the intrinsic physics of individual nodes of the network, but also what nodes talk to whom and how that communication influences nodal behaviors. The implication of such methodology is far-reaching -- one can begin to infer properties of networked systems with respect to network topology and/or external perturbations, answering hypotheticals such as ``what would happen if we remove this node?''
\end{quotation}

\section{Introduction} \label{sec:intro}

Networked dynamical systems, i.e., those described by the evolution of differential equations on a graph, are ubiquitous in science and engineering. When applied in the context of modeling dynamics of a population of units (e.g., a collection of fireflies or neural field activity), such systems exhibit a wide variety of population-wide complex behaviors including synchrony, pattern formation, chaos, and characteristic transitions between these states. 

Presently, we are concerned with the task of inferring (i) the physics of the individual units that make up the broader population, (ii) the underlying graphical structure of the units, and (iii) the coupling physics of a given networked dynamical system given nodal observations. Knowledge of these physics can allow for inference of system properties beyond future state predictions, including system-level stability and the emergence of coherent structures. Furthermore, in disambiguating between unit physics, coupling physics, and graphical interactions, one can project learned physics onto new, previously unexplored settings. For example, such tasks might include investigating system stability with respect to node deletions or other topological changes to the underlying graph.

In Section \ref{sec:background}, a brief overview of networked dynamical systems is presented in the context of state-of-the-art data-driven modeling efforts. In Section \ref{sec:methods}, the methodology for the present work is laid out. Numerical experiments are presented in Section \ref{sec:results} and is followed by discussion in Section \ref{sec:discussion}.

\begin{figure*}[]
        \centering
        \begin{overpic}[width=1.0\linewidth]{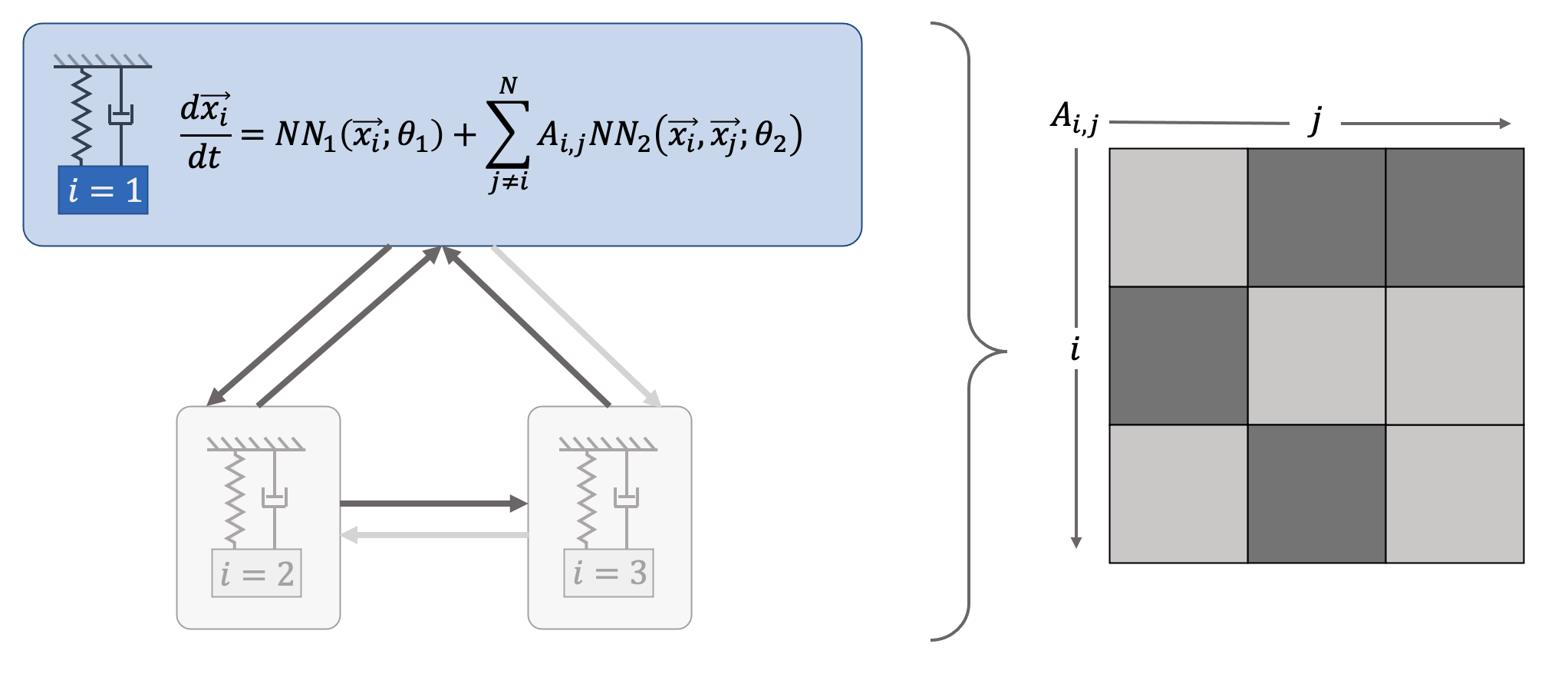}  
        \put(1,44){(a)}
        \put(67,44){(b)}
	    \end{overpic}  
	    \caption{In (a), a small network of unknown oscillators interact according to an adjacency matrix shown in (b). In this pictoral example, the nodes are assumed to be second-order oscillators with a directed graph structure, $A_{i,j}$. $\NN_1$ approximates the intrinsic nodal physics and $\NN_2$ approximates the coupling physics. The adjacency matrix in (b) is used in Section \ref{sec:inference}.}
		\label{fig:intro}
\end{figure*}

\section{Background} \label{sec:background}

Networked dynamical systems is an area of intense research in engineering and the sciences. Many physical systems feature a collection of individual units, or \textit{nodes}, possessing common intrinsic physics. The manner in which these nodes communicate - both \textit{how} nodes communicate and \textit{with whom} information is shared or directed - is of great interest, as these physics can inform system design and control. 

Examples in the physical sciences and engineering range from biological networks and pattern formation in neuroscience \cite{Earn1998,Klemm2005} to chemical reaction networks \cite{Ji_2021, Owoyele_2022} to thermoacoustics of combustors \cite{Dange2019,Premraj2021,Guan2022,Manoj2019}. Coupled dynamical systems have also been used successfully to model social and epidemiological networks \cite{Bakhtiarnia2020, Keeling2005}. Despite existing in disparate domains with distinct physics, these systems all can exhibit common qualitative behaviors like synchrony, pattern formation, and characteristic bifurcations between these states \cite{Xiao_2008, Strogatz_2001, Arenas_2008,Li2006}.

The Kuramoto family of coupled oscillators \cite{Kuramotoa} are the simplest nonlinear oscillator exhibiting many of these qualitative behaviors of physical systems. The generic Kuramoto oscillator possesses an assumed natural frequency $\omega$ and an all-to-all sinusoidal coupling:
\begin{equation}
    \frac{d\theta_i}{dt} = \omega + \frac{K}{N} \sum_{j=1}^{N}{\sin\left( \theta_i - \theta_j \right)} .
\end{equation}

The utility of the Kuramoto model is enshrined in its simplicity: basic modifications to the network topology (encoded through an \textit{adjacency} matrix) and coupling term can produce simple, interpretable models for chaos, turbulence, chimeras, and bifurcations between behaviors on an application-by-application basis. As a phase oscillator, the Kuramoto oscillator lacks the ability to capture higher-order dynamics. Nevertheless, the general prescription for a graph-coupled nonlinear dynamical system follows from that of the Kuramoto model:
\begin{equation}
        \frac{d\mathbf{x}_i}{dt} = \mathbf{f}(\mathbf{x}_i) + \sum_{j\neq i}^N \mathbf{A}_{i,j}\mathbf{g}\left(\mathbf{x}_i, \mathbf{x}_j\right) ,
\end{equation}
where $\mathbf{x}$ are the node states, $\mathbf{f}$ is the node physics, $\mathbf{g}$ is the coupling physics, and the adjacency matrix $\mathbf{A}_{i,j}$ directs the communication between nodes. 

\subsection{Related Work}
Many methods exist in data science and engineering to model nonlinear dynamical systems directly from observations. These methods range from so-called ``black-box'' approaches, where interpretability and/or physicality are not of concern, to physics-informed methods, where one uses prior knowledge of a system to construct, constrain, or otherwise regularize an optimization or regression task. Specifically for modeling dynamical systems, such methods include black-box neural networks for time-stepping (e.g., a multi-layer perceptron or recurrent neural network)~\cite{Skomski2021}, symbolic regression via a sparse selection of candidate functions \cite{Brunton_2016}, or Neural Ordinary Differential Equations (NODEs) \cite{Chen2018}. In the latter, unknown dynamical systems can be approximated with neural networks, mathematical terms known \textit{a priori} - potentially parameterized, or combinations of the two. Incorporating prior knowledge in NODEs -- e.g., a known term in the right-hand side (RHS) of the differential equation -- can bridge the gap between ``black-box'' modeling paradigms and from-first-principles forward modeling \cite{Rackauckas2020}. 

In the broader context of learning dynamics on graphs, Graph Neural Networks (GNNs) \cite{Scarselli2009,Poli2019,Rusch2022,Asikis_2022} aim to relate node features to those of their adjoining neighbors through diffusive processes (e.g. \textit{message passing} \cite{Duvenaud2015,Gilmer2017}) in an end-to-end differentiable algorithm. GNNs can process graph-structured data of arbitrary size and shape, making them very attractive for applications in physical chemistry and related fields. Although GNNs typically operate over static graphs, nodal hidden states are dynamic. These dynamics can be captured with traditional time-stepping models leading to natural extension of continuous graph neural networks (CGNNs)~\cite{CGNN2020} also referred to as graph neural ordinary differential equations (GDEs)~\cite{GDEs2019}.
Additionally, inference of a system's underlying graph can be accomplished within the GNN framework through generative models, i.e., graph variational autoencoders \cite{Kipf2016}.

The evolution of dynamics on graphs is intimately related to space-time dynamics governed by Partial Differential Equations (PDEs). Considering the discretization of space in a continuously varying spatial field, the states associated with the discrete nodes evolve according to some intrinsic physics subject to spatial interaction. In the case of a PDE, the spatial interaction follows the application (and truncation) of the Taylor series approximation for the differential operators on the discretized domain, resulting in sparse (and often symmetric) adjacency matrices. The corresponding data-driven method for approximating PDEs in a consistent manner is Neural PDEs -- a method-of-lines technique where the spatial interaction (stencil of the discretization) is learned through a convolutional neural network \cite{Rackauckas2020,Eliasof2021,Chamberlain2021}. 

Networked dynamical systems do not necessarily exist over an explicit spatial domain and instead have a sparse node-to-node coupling through fundamentally different physics (e.g., resistive coupling \cite{Wang2019}, or pulse coupling \cite{Abbott1993}). In such systems, learning the underlying graph is difficult. In \cite{Song2021}, a neural network is trained on a large number of samples ($10^6$) to output an adjacency matrix directly from observations of standard Kuramoto oscillators. Other network reconstruction tasks have been attempted with sparse regression techniques \cite{Han2015}, or by leveraging community structure \cite{Yan2012}. Compressive sensing techniques have been successfully used to identify and place missing nodes \cite{Su2012, Wu2012}. More recently, methods using the Gumbel Softmax technique \cite{Jang2016} have been used to propose adjacency matrices from data \cite{Zhang2019}, where jointly trained network proposal models and dynamics models work together to produce a candidate model for observed dynamics on a graph.

\subsection{Contributions}
In this work, we are motivated by modeling tasks in which network connectivity strongly influences qualitative behavior (i.e., limit cycling and bifurcation structure). To that end, we emphasize learning physically-consistent models that are interpretable and extensible. Our main contribution to this field is the construction and demonstration of a modeling framework in which one can recover system attributes required for informed modeling and control, namely (i) the intrinsic nodal physics, (ii) the coupling physics, and (iii) the underlying adjacency matrix. Our method utilizes universal differential equations to construct a unified approach for cooperatively learning the adjacency and physics components. We also show this method's generalization value by extending learned models to never-before-seen networks. 

\begin{figure}[]
        \centering
        \begin{overpic}[width=2.25in]{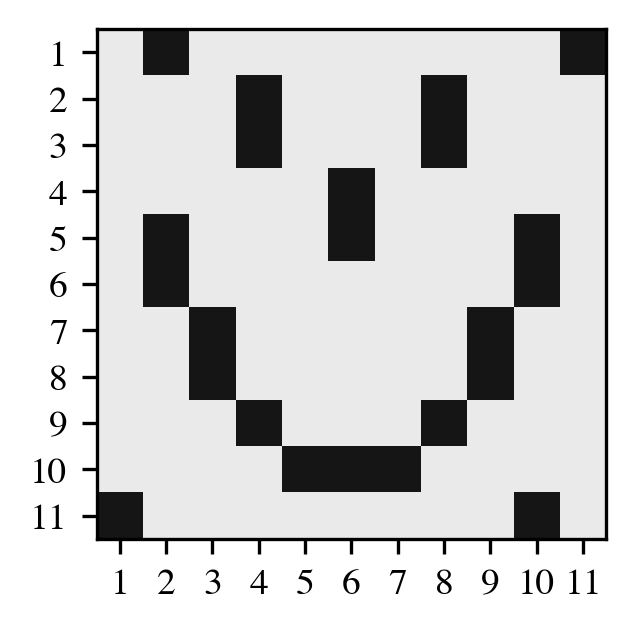}  
	    \end{overpic}  
	    \caption{Adjacency matrix for $11$-node directed network of nonlinear oscillators as defined in Eq. \ref{eq:system}. Dark elements correspond to $\mathbf{A}_{i,j} = 1$ with $\mathbf{A}_{i,j} = 0$ otherwise.}
		\label{fig:adjacency}
\end{figure}

\begin{figure}[]
        \centering
        \begin{overpic}[width=3in]{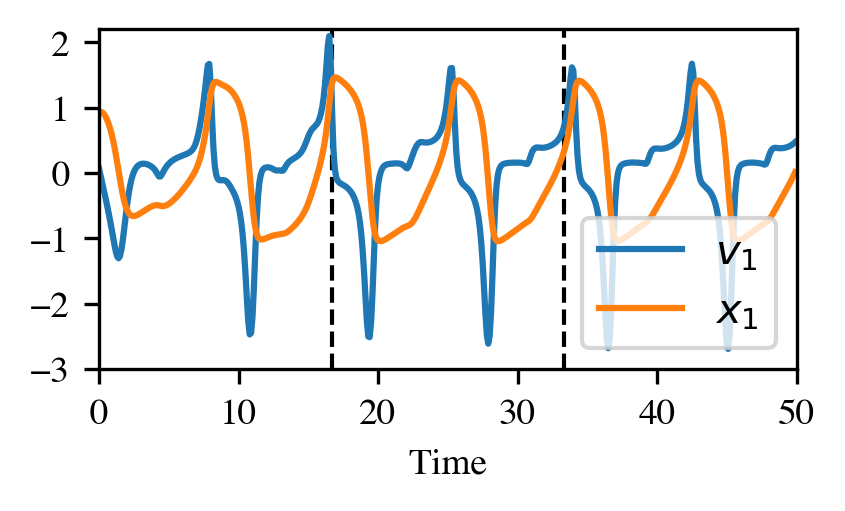}  
	    \end{overpic}  
	    \caption{Single node data ($i=1$) for $t \in [0,50)$ for position $x$ and velocity $v$. The vertical dashed lines indicate intervals for data splits into thirds for train-development-test sets.}
		\label{fig:single_node_data}
\end{figure}

\begin{figure*}[]
        \centering
        \begin{overpic}[width=6.69in]{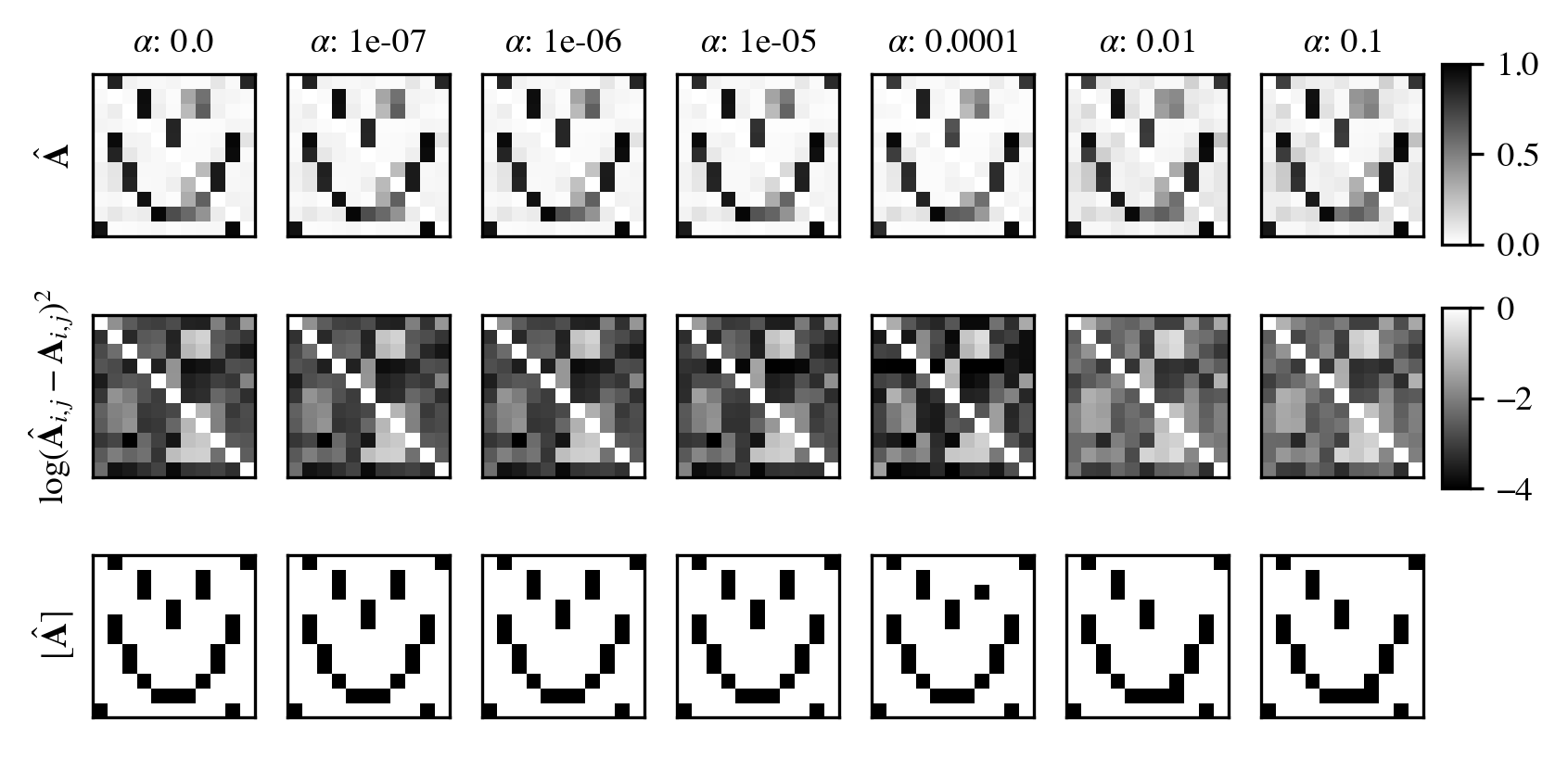}  
	    \end{overpic}  
	    \caption{Sweep of sparsity promotion penalty weight, $\alpha$, over several orders of magnitude. The top row is the visualization of the learned adjacency matrix ($\hat{\mathbf{A}}$) for the different values of $\alpha$. The middle row corresponds to the element-wise squared error between the actual and learned matrices. The last row shows the thresholding step, i.e., rounding to the nearest integer.}
		\label{fig:sweep}
\end{figure*}

\begin{figure}[]
        \centering
        \begin{overpic}[width=3in]{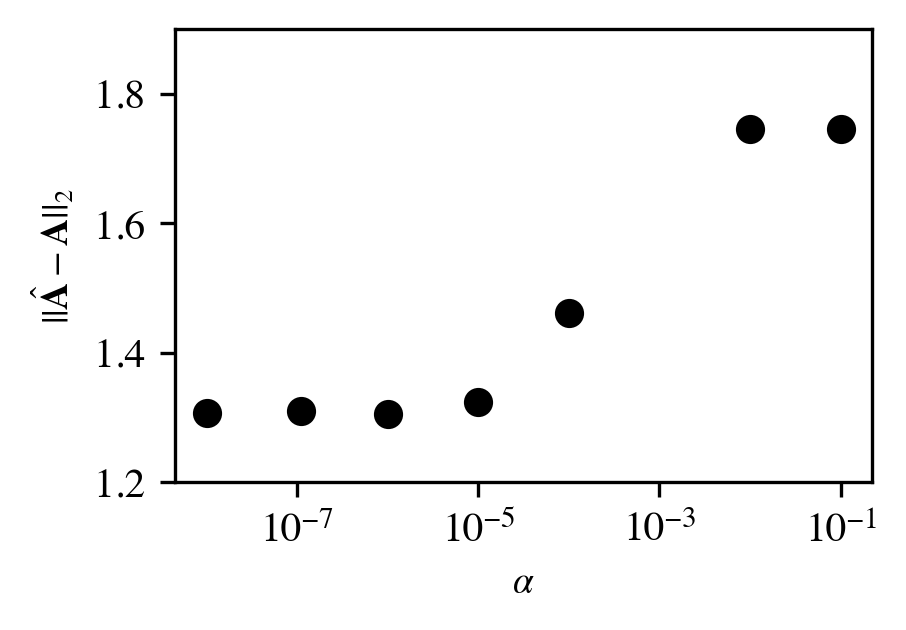}  
	    \end{overpic}  
	    \caption{L2 error between learned and actual adjacency matrix over sweep of sparsity promotion penalty weight, $\alpha$. The lowest error occurs at $\alpha = 1\cdot 10^{-6}$, although the error is insensitive to error until several orders of magnitude higher $\alpha$.}
		\label{fig:error_metric}
\end{figure}

\begin{figure*}[]
        \centering
        \begin{overpic}[width=6.69in]{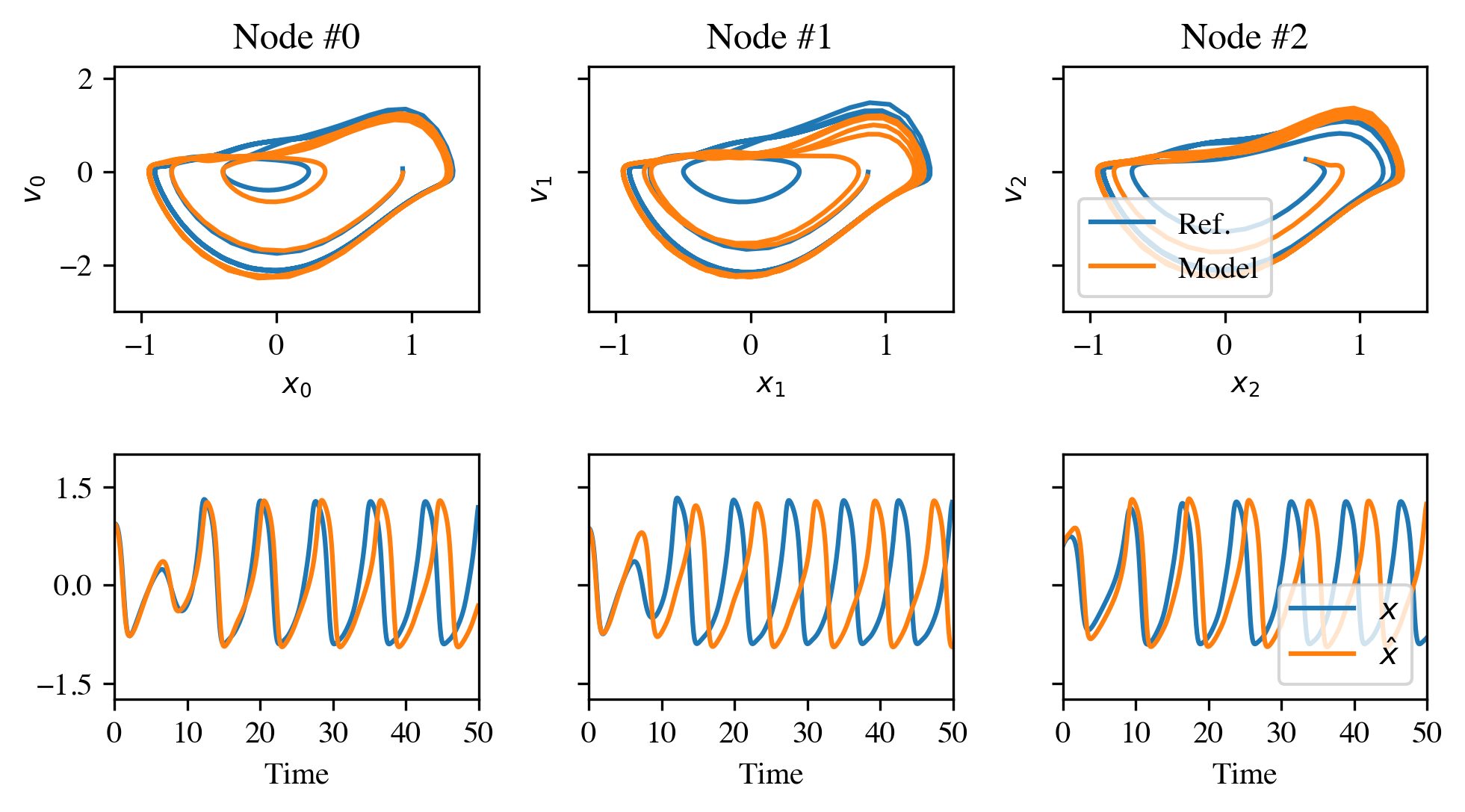}  
	    \end{overpic}  
	    \caption{The learned nodal and coupling physics from the coupled system with the 11-node directed network of Fig. \ref{eq:adjacency} is deployed on a never-before-seen 3-node network. The adjacency matrix for this network is shown in Fig. \ref{fig:intro}b. The learned physics models qualitatively capture the transients in the time signals, including the limit cycle behavior. }
		\label{fig:extrap}
\end{figure*}

\begin{table*}

\caption{Model selection metrics. Boldface represents the best (lowest) value among the candidate models.}

\label{tab:metrics}
\centering
\begin{tabular}{|cc|rrrrrrr|}
\hline
& & \multicolumn{7}{c|}{$\alpha$} \\

& & 0 & $1\cdot 10^{-7}$ &  $1\cdot 10^{-6}$ &  $1\cdot 10^{-5}$ &  $1\cdot 10^{-4}$ &  $1\cdot 10^{-2}$ &  $1\cdot 10^{-1}$ \\ 

\hline
\multicolumn{2}{|c|}{${||A||}_1$} & 24.59 & 24.58 & 24.46 & 23.82 & 22.16 & 26.74 & 26.74 \\ 

\hline
\multirow{3}{*}{\rotatebox[origin=c]{90}{\makecell{$N_f$-step \\MSE}}} & Train &  0.00789 &  0.01108 &  0.00761 &  \textbf{0.00643} &  0.00837 &  0.01392 &  0.01392 \\
& Dev. &  0.01019 &  0.01442 &  0.00968 &  \textbf{0.00911} &  0.01153 &  0.01790 &  0.01790 \\
& Test &  0.01002 &  0.01415 &  0.00951 &  \textbf{0.00905} &  0.01148 &  0.01757 &  0.01757 \\ 

\hline
%\multirow{3}{*}{\rotatebox[origin=c]{90}{\makecell{Loop\\ MSE}}} & Train & 1.50524 &  1.61105 &  1.52188 &  0.76406 &  0.71677 &  1.78833 &  1.78833 \\
%& Dev. &  1.56974 &  1.67695 &  1.59397 &  1.03009 &  0.91528 &  1.91965 &  1.91965 \\
%& Test &  1.53878 &  1.61713 &  1.56156 &  0.99795 &  0.89865 &  1.83894 &  1.83894 \\
%\hline
\end{tabular}
\end{table*}

\section{Methodology} \label{sec:methods}
Our goal is to identify a networked nonlinear dynamical system from measurements of nodal states over time. We begin by defining $N$ assumed-identical nodes with states $\mathbf{x}_i(t) \in \mathbb{R}^d$ where $i = (1, 2, ..., N)$, $d$ is the state dimension of the $i$-th node, and $N\cdot d$ is the system state dimension. The connectivity of the nodes is prescribed by an assumed time-invariant adjacency matrix $\mathbf{A} \in \mathbb{R}^{N  \times  N}$, where entry $\mathbf{A}_{i,j}$ prescribes the learnable (parameterized by $\theta_g$) interaction $\mathbf{g}(\mathbf{x}_i, \mathbf{x}_j;\theta_g): \mathbb{R}^{2d} \rightarrow \mathbb{R}^d$ between nodes $i$ and $j$. This interaction is assumed to be directed ($\mathbf{A}$ not necessarily symmetric), unweighted (meaning our final adjacency during deployment should have binary elements), and sparse.

Thus, we define an expected dynamical system (Figure \ref{fig:intro}) of the form:
\begin{equation} \label{eq:node}
        \frac{d\mathbf{x}_i}{dt} = \mathbf{f}\left(\mathbf{x}_i;\theta_f\right) + \sum_{j\neq i}^N \mathbf{A}_{i,j}\mathbf{g}\left(\mathbf{x}_i, \mathbf{x}_j;\theta_g\right) = \mathbf{F}\left(\mathbf{x};\Theta\right),
\end{equation}
where $\mathbf{f}(\mathbf{x}_i,\theta_f) : \mathbb{R}^d \rightarrow \mathbb{R}^d$ is the learnable intrinsic physics of the nodes parameterized by $\theta_f$ and $\Theta$ is the concatenated vector of model parameters. In general, $\mathbf{f}$ and $\mathbf{g}$ can be any differentiable surrogate model (e.g. a neural network) mapping to the full nodal state dimension $d$. However, these can be constructed to satisfy \textit{a priori} knowledge of a system; e.g. a parameterized equation, composite function, or other closure-type model \cite{Rackauckas2020}.

A Neural ODE evolve states $\mathbf{x}$ through time $t$ according to the transition dynamics specified by the RHS of Eq. \ref{eq:node}. The evolution of $\mathbf{x}$ corresponds to the initial value problem:
\begin{equation} \label{eq:ivp}
    \mathbf{x}\left(t_{end}\right) = \mathbf{x}\left(t_0\right) + \int_{t_0}^{t_{end}} \mathbf{F}\left( \mathbf{x};\Theta \right) dt .
\end{equation}
Equation \ref{eq:ivp} can be solved with standard numerical ODE solvers,
\begin{equation} \label{eq:odesolve}
    \mathbf{x}\left(t_{end}\right) = \ODESolve \left( \mathbf{F}\left( \mathbf{x};\Theta \right), \mathbf{x}_0, t_0, t_{end} \right) .
\end{equation}
Equation \ref{eq:odesolve} can be made fully differentiable one of two ways: (i) via the adjoint sensitivity method, whereby a second differential equation is formulated and solved for the evolution of an error metric backward in time from $t_{end}$ to $t_0$ (see Chen et al. \cite{Chen2018}), or (ii) through the computational graph of the forward pass of the ODE solver (e.g. composition of differentiable elementary operations of the fourth-order Runge-Kutta integrator, for example). The latter method is employed in this work for its ease of implementation.

A necessary output of the training routine is a adjacency matrix with connections between nodes chosen to fit the data. Assumed is that the network to be modeled is unweighted, i.e., nodes $i$ and $j$ are either ``connected'' with $A_{i,j} = 1$ or ``disconnected'' with  $A_{i,j} = 0$. Several methods exist for making such a selection task differentiable (e.g. the Gumbel-Softmax trick). In this work, to begin training, we instead assume all-to-all coupling activated with a sigmoid:
\begin{equation} \label{eq:adjacency}
    \mathbf{A} = \sigma \left( \left[ 0 \right] - I/\epsilon \right),
\end{equation}
with small number $\epsilon$, such that the off-diagonal entries of $\mathbf{A}$ default to a value of one-half and are bounded to $\mathbf{A}_{i,j} \in (0,1)$. After training is complete, the entries of $\mathbf{A}$ are thresholded to either $0$ or $1$ based on some criteria; for the present study, this is simply a rounding operation on the elements of $\mathbf{A}$ post-training. 

The training data is formatted in a batched tuple with $N_f$-step moving-window forecasting horizon 
\begin{equation}
    \mathcal{D}=\{ \mathbf{X}^{(k)} (t), \mathbf{X}^{(k)} (t+\Delta), ..., \mathbf{X}^{(k)} (t+N_f \Delta)\},
\end{equation}
where $\mathbf{X}=\{\mathbf{x}_1,...,\mathbf{x}_N\}$; $k$ is the $k$-th batch ($k=1,...,N_b$); $\Delta$ is the size of discrete timesteps. During training, the surrogate model with parameters $\boldsymbol{\Theta}$ is asked to predict the following $N_f-1$ timesteps provided $\mathbf{X}^{(k)} (t)$ for each batch $k$. Therefore, the loss function to be minimized could be constructed as the mean squared error between predictions and measurements, given as
\begin{equation} \label{eq:loss}
    \mathcal{L} (\boldsymbol{\Theta}, \mathbf{A}) = \frac{1}{N_b \cdot N_f \cdot N \cdot d} || \hat{\mathbf{X}}^{(k)} (t; \boldsymbol{\Theta}) - \mathbf{X}^{(k)} (t)||^2_2 + \alpha|| \mathbf{A}||_1, 
\end{equation}
where $\hat{\mathbf{X}}^{(k)} (t; \boldsymbol{\Theta})$ is prediction for all system states at time $t$ in batch $k$; $|| \mathbf{A}||_1$ is a $\ell_1$ penalty intended to promote sparse connections in the adjacency matrix; and $\alpha$ is the weighting coefficient for this $\ell_1$ penalty.

\section{Numerical Experiments} \label{sec:results}

Presented are two numerical experiments designed to show the utility of the presented methods. First, in Section \ref{sec:modeling}, we perform the data-driven modeling of a coupled oscillator system with unknown adjacency, node physics, and coupling physics. Once these physics are learned, we shift our attention to inferring system behavior on a never-before-seen network in Section \ref{sec:inference}. 

\subsection{Modeling Coupled Nonlinear Oscillators} \label{sec:modeling}
Consider a canonical nonlinear oscillator:
\begin{align}
    \frac{d^2 x_i}{dt^2} = & -x_i - a_1 \frac{dx_i}{dt} \left( a_2 x_i^4 - a_3 x_i + a_4 \right) \\
     & + \sum_{j=1}^{N} \mathbf{A}_{i,j} \left( \frac{d^2 x_i}{dt^2} - \frac{d^2 x_j}{dt^2} \right) \nonumber,
\end{align}
or, written as a system of first-order ordinary differential equations:
\begin{align} \label{eq:system}
    \frac{d x_i}{dt} = & v_i \\
    \frac{d v_i}{dt} = & -x_i - a_1 v_i \left( a_2 x_i^4 - a_3 x_i + a_4 \right) \nonumber\\
     & + \sum_{j=1}^{N} \mathbf{A}_{i,j} \left( v_i - v_j \right) \nonumber,
\end{align}
with $N=11$, $a_{1,...,4} = \{0.2, 11.0, 11.0, 1.0\}$, $\mu=0.2$, and adjacency matrix specified by that of Fig. \ref{fig:adjacency}. 

Depending on parameter values ($a_{1,...,4}$ and $\mu$), the system detailed in Eq. \ref{eq:system} exhibits unique behaviors; for the specific parameter values and network detailed, this system shows limit cycle behavior, as seen in the data generated for this study and shown in Fig. \ref{fig:single_node_data}. For demonstration of our methodology, the system in Eq.  \ref{eq:system} was simulated from a random initial condition for 500 timesteps with temporal resolution of $\Delta t = 0.1$. This constitutes our dataset for this task.

We aim to fit the notional coupled dynamical system (Eq. \ref{eq:node}) to this dataset. The trainable surrogate for nodal intrinsic physics is $\NN_1(x_i, v_i; \theta_1): \mathbb{R}^2 \rightarrow \mathbb{R}^1$, and similarly for the interaction physics: $\NN_2(\mathbf{x}_i, \mathbf{x}_j; \theta_2): \mathbb{R}^4 \rightarrow \mathbb{R}^1$. Thus, by assuming second-order physics, our model takes the form:
\begin{align} \label{eq:model}
    \frac{d x_i}{dt} = & v_i \\
    \frac{d v_i}{dt} = & \NN_1(x_i, v_i; \theta_1) + \sum_{j \neq i}^{N} A_{i,j} \NN_2(\mathbf{x}_i, \mathbf{x}_j; \theta_2) \nonumber \\
    A_{i,j}(\theta_3) = & \sigma \left( \vectorize^{-1}(\theta_3) - I/\epsilon \right) \nonumber ,
\end{align}
where the $\vectorize^{-1} \left( \cdot \right)$ operator reshapes the vector $\theta_3 \in \mathbb{R}^{N^2}$ to $\mathbb{R^{N\times N}}$. The network $\NN_1$ is chosen to be a feed-forward network with two hidden layers of 50 nodes each with Leaky ReLU activation functions. $\NN_2$ has a single hidden layer of size 4 with Leaky ReLU activation applied. 

The integration scheme for solving Eq. \ref{eq:model} is a straightforward fourth-order Runge-Kutta integrator. The forecast horizon is set to $N_f = 5$ during training. The optimizer used is Adam with a fixed learning rate of 0.02. First, the model parameters are optimized for 1000 epochs to minimize the loss as specified by Eq. \ref{eq:loss} with $\alpha = 0$; that is, without penalizing adjacency sparsity. After 1000 epochs, the sparsity penalty weight $\alpha$ is increased to a fixed value and the optimization continues with another 1000 epochs. 

Table \ref{tab:metrics} summarizes the error metrics across the training, development, and test data splits (Figure \ref{fig:single_node_data}) through a sweep of values for $\alpha$ (log-10 scale). For the purposes of this work, the model selection criteria is the trained model with the lowest Mean-Squared-Error (MSE) on the development data split. For this set of experiments, a value of $\alpha = 1 \cdot 10^{-5}$ resulted in the lowest error for this metric. As shown in Fig. \ref{fig:error_metric}, the corresponding learned adjacency matrix identically matches the ground-truth adjacency matrix after a thresholding operation. In general, as the sparsity promotion penalty increases, several elements in the adjacency matrix do not meet the threshold and are set to zero. Note, however, that the error between the learned and true adjacency matrices is relatively insensitive to $\alpha$ for low values of $\alpha$ (less than $\alpha = 1 \cdot 10^{-5}$) for our chosen training procedure.

\subsection{Inference} \label{sec:inference}
The trained model of Section \ref{sec:results} contains neural network approximations for the intrinsic nodal physics ($\NN_1$) and nodal coupling physics ($\NN_2$). By virtue of the construction of our model as a system of ODEs, these approximations should hold regardless of the underlying graph structure of the system. To test this, these physics approximations are deployed on a small, never-before-seen three-node network whose adjacency matrix is given by Fig. \ref{fig:intro}b. 

Figure \ref{fig:extrap} shows phase portraits and time series for the three nodes of the new network. In the top row, position $x_i$ and velocity $v_i$ show the clear limit cycle behavior for both a ground-truth simulation and the roll-out of the learned model. In the second row, the time series reconstruction (displayed is $x_i$) from a random initial condition for each node is shown. Although this network was not included in the training data, the learned nodal and coupling physics approximations perform well in reproducing the major qualitative behavior of the system, including limit cycling and following transient behavior of the nodal states as they approach the limit cycle (e.g. for time $0<t<20$ for each of the nodes). 

\section{Discussion} \label{sec:discussion}
Although an exhaustive set of numerical experiments was not presented, we have shown the general application and utility of universal differential equations for structural inference of networked dynamical systems. Specifically, we have shown that one can approximate (i) nodal physics, (ii) coupling physics, and (iii) the underlying graphical structure of such systems directly from data. Furthermore, we have demonstrated that these learned physics can be successfully applied to never-before-seen networks and capture the general qualitative behavior of these systems.

Several modeling choices were made to deliver the results contained within this work; examples include the neural network architecture(s), the network characteristics (weighted vs. unweighted, directed vs. undirected, etc.), and the neural ODE model structure and integration scheme, among may others. While each of these aspects of the modeling process can have large impacts on model quality and performance, in this section detailed are those deserving of extra attention. Specifically we focus on (i) potential shortcomings of the framework and avoidance strategies and (ii) applications and future development work.

\subsection*{Shortcomings}
The model form of Eq. \ref{eq:model} is generic in that assumed is second-order dynamics with additive coupling - no other structure or physical prior is enforced in the model construction. This contrasts so-called ``gray'' or ``white''-box models in which known terms are inserted into the model as \textit{a priori} knowledge. In this context, such prior knowledge can take the form of known nodal physics, coupling physics, and/or adjacency matrix. The class of models presented in this work fall into the category of ``black''-box models, which are prone to overfitting, especially in the context of neural ODEs. Model performance will likely suffer for systems with (i) weakly nonlinear dynamics (e.g. a Kuramoto oscillator network with small coupling coefficient) or (ii) data existing solely on the system's attractor. In the case of the former, disambiguation of distinct signals from time series data may be difficult or impossible, especially in the situation of global synchrony (a population of homogeneous Kuramoto oscillators). Likewise, should observations of a system exist exclusively on that system's attractor, it is likely that a neural ODE-based model may learn exclusively the shape of the attractor in phase space - learned models are unlikely to generalize well off of the attractor. Inclusion of prior knowledge will strongly regularize the neural ODE such that learned closure terms or parameters will generalize better than their black-box counterparts. For example, reformulating the learning task to learn exclusively the adjacency matrix of a networked dynamical system given known nodal and coupling physics will generally produce a better model than the black-box counterpart.

Similarly, the networks presented in this work have $N=11$ for Section \ref{sec:modeling} and $N=3$ for \ref{sec:inference} with uniform nodal physics. These are small networks with sparse connections, though for the presented training routine, assumed is all-to-all coupling meaning that the number of trainable parameters (for $\mathbf{A}$) exhibits quadratic growth with the number of nodes, underscoring the need for regularization during training (e.g. $\ell_1$ penalty for $\mathbf{A}$).

\subsection*{Applications and Future Work}
This work was motivated by the open problem of inferring the graphical structure of networked dynamical systems from data. Following the layout of the presented results, this modeling framework can be directly applied to (i) separate nodal and coupling physics from nodal interactions and (ii) estimate the stability properties of an observed system with respect to perturbations in the underlying network. 

The class of physical systems amenable to the presented techniques is wide-ranging. In closed combustion systems, thermoacoustic instabilities can arise from coupling of injectors with fluid dynamic processes of the combustor-injector placement and control is an area of intense study, which amounts to perturbations of an adjacency matrix of the coupled injectors\cite{Dange2019,Manoj2019,Koch2020}. In power system network dynamics, nodes obey the so-called swing equations for synchronous operation of equipment. Such physics can be approximated by a network of coupled oscillators, from which one can estimate stability properties of the system directly from data. Extending the method further, time-varying adjacency matrices can be incorporated to model molecular dynamics, whereby the underlying graphical structure evolves with the nodes. In general, broad applicability can be had by encoding physical priors. For example, one could specify a known energy input/output balance in a model for a dissipative system (e.g. building thermal dynamics and other lumped capacitance-type models) such that the model can be readily constructed and analyzed in the context of a particular graph structure. Our modeling framework can inform system design and control in these contexts from a purely data-driven perspective, e.g., the modification or evolution of the graph structure and system control, among other tasks.

\section{Conclusion} \label{sec:conclusion}
We have demonstrated the use of universal differential equations for the data-driven modeling of networked dynamical systems. Specifically, we have demonstrated the capability to extract (i) intrinsic nodal physics, (ii) coupling physics, and (iii) the underlying graphical structure of the system. This is accomplished by assuming a generic differential equation form that can be applied to a wide variety of networked dynamical systems. The modeling framework was demonstrated on an 11-node-directed network of nonlinear oscillators. The utility of this framework lies in the disambiguation of the physics from network structure. In doing so, the physics can be applied to a network of arbitrary topology, allowing one to investigate system properties with respect to perturbations in the underlying network.

\section*{Acknowledgements}
This work was supported through the Data Model Converge (DMC) initiative at the Pacific Northwest National Laboratory (PNNL). PNNL is a multi-program national laboratory operated for the U.S. Department of Energy (DOE) by Battelle Memorial Institute under Contract No. DE-AC05-76RL0-1830.

\section*{Conflicts of Interest}
The authors have no conflicts to disclose.

\section*{Author's Contributions}
\textbf{J. Koch:} conceptualization, software, investigation, visualization, draft preparation. \textbf{Z. Chen:} software, investigation, visualization, draft preparation. \textbf{A. Tuor:} software, supervision. \textbf{J. Drgona:} software, draft editing, supervision. \textbf{D. Vrabie:} funding acquisition, advisor.

\section*{Data Availability}
Data available from the authors upon reasonable request.

\section*{References}
\bibliography{node_bib}% Produces the bibliography via BibTeX.

\end{document}